# CONSTRUCTION OF COCYCLE BICROSSPRODUCTS BY TWISTING

**G. BOGDANOFF, I. BOGDANOFF**
*Megatrend University of Applied Sciences*
*Belgrade - Serbia*
igor.bogdanoff@megatrend.edu.rs
grichka.bogdanoff@megatrend.edu.rs

**Abstract :** We consider a generalisation of the Majid's mirror product $M(H) = H^{op} \triangleright\blacktriangleleft H$ of a Hopf algebra $H$, when one of the components of the product is replaced by a twist. This leads to a new "twisted mirror product" construction for cocycle bicrossproduct Hopf algebra of the form $M_\chi(H) = H^{op} {}^\psi\triangleright\blacktriangleleft H_\chi$, where $\chi$ is a 2-cocycle of "twisting" type. As an application, we obtain a canonical bicrossproduct $M(H) = H^\psi \triangleright\blacktriangleleft H$ for any Hopf algebra $H$.

## 1. INTRODUCTION

Let $H$ be a Hopf algebra. In [1,2] was introduced a general construction of bicrossproduct Hopf algebras of the form $M(H) = H^{op} \triangleright\blacktriangleleft H$ where the algebra is a semidirect product by a version of the quantum adjoint action of $H^{op}$ on $H$ and the coalgebra is a semidirect coproduct by an adjoint coaction of $H$ on $H^{op}$. Although isomorphic to a tensor product, this "mirror product" $\triangleright$ [3] provided on the one hand examples of general non Abelian type (where $H$ is neither commutative nor cocommutative) and on the other hand it corresponds under a general semidualisaiton construction [3] (in the finite dimensional case) to the Drinfeld quantum double in the double cross product form $H^{*op} \triangleright\triangleleft H$.

In this note, we point out a useful generalisation of this mirror construction, originally motivated from physics [4,5,6]. Namely, if $\chi$ is a 2-cycle in $H \otimes H$ in the sense of Drinfeld [7], then it is known that there is a new Hopf algebra $H\chi$ where the coproduct is conjugated by $\chi$ (see [8] for the explicitly Hopf algebra case). We show in section 2 that this cocycle twist can also be used to build a dual cocycle in the sense of cocycle bicrossproducts [2] leading to a "twisted mirror product" and cocycle bicrossproduct extension

$$H\chi \to M\chi(H) \to H^{op} \ ; \quad M\chi(H) = H^{op\,\psi} \triangleright\!\!\blacktriangleleft H\chi \qquad (1)$$

Here $H\chi$ appears as a subHopf algebra.

In section 3, we obtain a related construction $\overline{M}(H) = H^{\psi} \triangleright\!\!\blacktriangleleft H$. If $H$ were quasitriangular in the sense of Drinfeld [9], then this would be a special case of the above with $\chi$ built from the quasitriangular structure. However, the construction is more general. Basically, the introduction of a certain cocycle $\psi$ allows us to replace $H^{op}$ by $H$ in Majid's construction.

Section 4 contains some concluding remarks. Let us also mention the physical motivation [4,5,6], though this will not be relevant here. Namely, it is known that the q-Lorentz group $Uq(so(3, 1))$ is a twisting of the q-Euclidean group $Uq(so(4))$ (see [10]). Hence, the above $M\chi(Uq(so(3, 1)))$ provides a single Hopf algebra of the form $Uq(so(4))^{op\,\psi} \triangleright\!\!\blacktriangleleft Uq(so(3, 1))$ unifying the Lorentzian and Euclidean ones. This result provides a possibly interesting algebraic basis for the hypothesis of fluctuations of the signature of the spacetime metric at the Planck scale proposed in (4).

## 2. THE COCYCLE BICROSSPRODUCT $M\chi(H) = H^{op\,\psi} \triangleright\!\!\blacktriangleleft H\chi$

We recall that the "mirror product" $M(H) = H^{op} \triangleright\!\!\blacktriangleleft H$ is based on the actions and coactions :

$$a \triangleleft h = h_{(1)}\, a\, Sh_{(2)}, \quad \forall\ a \in H,\ h \in H^{op}$$

$$\beta(h) = h_{(1)}\, Sh_{(3)} \otimes h_{(2)}, \quad \forall\ h \in H^{op} \qquad (2)$$

We use here the Sweedler notation (preferably refer Sweedlers book on Hopf algebras) for the coproduct, with summation of terms understood. We use conventions for Hopf algebras and their bicrossproducts as in [11].

One then makes the semidirect product and coproduct by these in the usual way for Hopf algebras. We now prove a generalisation of this mirror product, when one of the components of the product is replaced by a twist. Then we get :

**Theorem 2.1** *Let $H$ be a Hopf algebra equiped with a bijective antipode, $\chi \in H \otimes H$ be a 2-cocycle and $H_\chi$ the twisted Hopf algebra of Drinfeld associated with the Hopf algebra H. Then there exists a cocycle bicrossproduct of the form*

$$M_\chi(H) = H^{op\,\psi} {\triangleright\!\!\blacktriangleleft}\, H_\chi$$

*where $\chi = \chi_{(1)} \otimes \chi_{(2)}$ is an explicit notation (summation understood), and similarly for $\chi_{\{-1\}}$. We have an extension of Hopf algebras*

$$H_\chi \to M_\chi(H) \to H^{op}$$

*Furthermore, $M_\chi(H) = H^{op} \otimes H_\chi$ as a Hopf algebra by $h_{(1)} \otimes h_{(2)} a \leftarrowtail h \otimes a$*

**Proof** We verify the conditions for the cocycle bicrossproduct $H^{op\,\psi} {\triangleright\!\!\blacktriangleleft}\, H_\chi$ with $\psi$ set out as above and theorem 6.3.9 in [11]. We note that $H^{op}$ plays the role of $H$ there, and $H_\chi$ plays the role of $A$ there. Note also that $H_\chi$ has the same algebra as $H$ and hence remains an $H^{op}$ - module algebra, as in the usual construction $M(H) = H^{op} {\triangleright\!\!\blacktriangleleft}\, H$. Next, we have:

$$(id \otimes \beta) \circ \beta(h_{(1)}) \psi_{12}(h_{(2)}) = (h_{(1)(1)} Sh_{(1)(3)} \otimes h_{(1)(2)(1)} Sh_{(1)(2)(3)} \otimes h_{(1)(2)(2)}) \psi_{12}(h_{(2)})$$

$$= h_{(1)(1)} Sh_{(1)(3)} h_{(2)} \chi^{(1)} Sh_{(5)} \chi^{-(1)} \otimes h_{(1)(2)(1)} Sh_{(1)(2)(3)} h_{(3)} \chi^{(2)} Sh_{(4)} \chi^{-(2)} \otimes h_{(1)(2)(2)}$$

$$= h_{(1)} Sh_{(5)} h_{(6)} \chi^{(1)} Sh_{(9)} \chi^{-(1)} \otimes h_{(2)} Sh_{(4)} h_{(7)} \chi^{(2)} Sh_{(8)} \chi^{-(2)} \otimes h_{(3)}$$

$$= h_{(1)} \chi^{(1)} Sh_{(5)} \chi^{-(1)} \otimes h_{(2)} \chi^{(2)} Sh_{(4)} \chi^{-(2)} \otimes h_{(3)}$$

whereas

$$\psi_{12}(h_{(1)})(\Delta_\chi \otimes id)\beta(h_{(2)}) = \psi_{12}(h_{(1)})(\chi^{(1)}h_{(2)(1)(1)}Sh_{(2)(3)(2)}\chi^{-(1)} \otimes \chi^{(2)}h_{(2)(1)(2)}Sh_{(2)(3)(1)}\chi^{-(2)} \otimes h_{(2)(2)})$$

$$= h_{(1)(1)}\chi^{(1)}Sh_{(1)(4)}h_{(2)(1)(1)}Sh_{(2)(3)(2)}\chi^{-(1)} \otimes h_{(1)(2)}\chi^{(2)}Sh_{(1)(3)}h_{(2)(1)(2)}Sh_{(2)(3)(1)}\chi^{-(2)} \otimes h_{(2)(2)}$$

$$= h_{(1)(1)}\chi^{(1)}Sh_{(1)(4)}h_{(2)(1)(1)}Sh_{(2)(3)(2)}\chi^{-(1)} \otimes h_{(1)(2)}\chi^{(2)}Sh_{(1)(3)}h_{(2)(1)(2)}Sh_{(2)(3)(1)}\chi^{-(2)}h_{(2)(2)}$$

$$= h_{(1)}\chi^{(1)}Sh_{(4)}h_{(5)}Sh_{(9)}\chi^{-(1)} \otimes h_{(2)}\chi^{(2)}Sh_{(3)}h_{(6)}Sh_{(8)}\chi^{-(2)} \otimes h_{(7)}$$

$$= h_{(1)}\chi^{(1)}Sh_{(5)}\chi^{-(1)} \otimes h_{(2)}\chi^{(2)}Sh_{(4)}\chi^{-(2)} \otimes h_{(3)}$$

as required. We have used here the elementary properties of Hopf algebras. Also, using the cocycle property of $\chi$, we can see that $\psi$ is a cocycle in the required sense, such that $H^{op}$ becomes a cocycle $H_\chi$ - comodule. It is clear immediately that the resulting cocycle coaction respects the coproduct of $H^{op}$ in as much as these maps are the same as for $M(H) = H^{op} \bowtie H$. Thus :

$$(id \otimes \psi) \circ \beta(h_{(1)})((id \otimes \Delta_\chi)\psi(h_{(2)})) = (h_{(1)(1)}Sh_{(1)(3)} \otimes \psi(h_{(1)(2)}))\chi_{23}((id \otimes \Delta)\psi(h_{(2)}))\chi_{23}^{-1}$$

$$= h_{(1)(1)}Sh_{(1)(3)}h_{(2)(1)}\chi^{(1)'}Sh_{(2)(4)}\chi^{-(1)'} \otimes h_{(1)(2)(1)}\chi^{(1)}Sh_{(1)(2)(4)}(h_{(2)(2)}\chi^{(2)'}Sh_{(2)(3)}\chi^{-(2)'})_{(1)}\chi^{-(1)} \otimes$$
$$h_{(1)(2)(2)}\chi^{(2)}Sh_{(1)(2)(3)}(h_{(2)(2)}\chi^{(2)'}Sh_{(2)(3)}\chi^{-(2)'})_{(2)}\chi^{-(2)}$$

$$= h_{(1)}Sh_{(6)}h_{(7)}\chi^{(1)'}Sh_{(12)}\chi^{-(1)'} \otimes h_{(2)}\chi^{(1)}Sh_{(5)}h_{(8)}\chi^{(2)'}{}_{(1)}Sh_{(11)}\chi^{-(2)'}{}_{(1)}\chi^{-(1)} \otimes$$
$$h_{(3)}\chi^{(2)}Sh_{(4)}h_{(9)}\chi^{(2)'}{}_{(2)}Sh_{(10)}\chi^{-(2)'}{}_{(2)}\chi^{-(2)}$$

$$= h_{(1)}\chi^{(1)'}Sh_{(6)}\chi^{-(1)'} \otimes h_{(2)}\chi^{(1)}\chi^{(2)'}{}_{(1)}Sh_{(5)}\chi^{-(2)'}{}_{(1)}\chi^{-(1)} \otimes h_{(3)}\chi^{(2)}\chi^{(2)'}{}_{(2)}Sh_{(4)}\chi^{-(2)'}{}_{(2)}\chi^{-(2)}$$

where $\chi'$ is another copy of $\chi$, whereas

$$\psi_{12}(h_{(1)})(\Delta_\chi \otimes id)\psi(h_{(2)}) = \psi_{12}(h_{(1)})\chi_{12}((\Delta \otimes id)\psi(h_{(2)}))\chi_{12}^{-1}$$

$$= h_{(1)(1)}\chi^{(1)}Sh_{(1)(4)}(h_{(2)(1)}\chi^{(1)'}Sh_{(2)(4)}\chi^{-(1)'})_{(1)}\chi^{-(1)} \otimes$$
$$h_{(1)(2)}\chi^{(2)}Sh_{(1)(3)}(h_{(2)(1)}\chi^{(1)'}Sh_{(2)(4)}\chi^{-(1)'})_{(2)}\chi^{-(2)} \otimes h_{(2)(2)}\chi^{(2)'}Sh_{(2)(3)}\chi^{-(2)'}$$

$$= h_{(1)}\chi^{(1)}Sh_{(4)}h_{(5)}\chi^{(1)'}_{(1)}Sh_{(10)}\chi^{-(1)'}_{(1)}\chi^{-(1)} \otimes$$
$$h_{(2)}\chi^{(2)}Sh_{(3)}h_{(6)}\chi^{(1)'}_{(2)}Sh_{(9)}\chi^{-(1)'}_{(2)}\chi^{-(2)} \otimes h_{(7)}\chi^{(2)'}Sh_{(8)}\chi^{-(2)'}$$

$$= h_{(1)}\chi^{(1)}\chi^{(1)'}_{(1)}Sh_{(6)}\chi^{-(1)'}_{(1)}\chi^{-(1)} \otimes h_{(2)}\chi^{(2)}\chi^{(1)'}_{(2)}Sh_{(5)}\chi^{-(1)'}_{(2)}\chi^{-(2)} \otimes h_{(3)}\chi^{(2)'}Sh_{(4)}\chi^{-(2)'}$$

which is equal to the above expression because of the cocycle axiom for $\chi$ and its corresponding version for $\chi^{-1}$. We thus obtain a cocycle coalgebraic cross coproduct of the form $H^{\text{op}} \ {}^{\psi}\!\!\triangleright\!\blacktriangleleft H_\chi$ and a crossproduct $H^{\text{op}} \triangleright\!\!< H_\chi$. We can similarly verify all of the required compatibility conditions in theorem 6.3.9 of [10] and hence have a bicrossproduct Hopf algebra. Alternatively, we note that we have an algebra isomorphism

$$\theta : H^{\text{op}} \otimes H_\chi \xrightarrow{\sim} H^{\text{op}} \triangleright\!\blacktriangleleft H_\chi \qquad (3)$$

$$\theta(h \otimes a) = h_{(1)} \otimes h_{(2)}a \qquad (4)$$

because the algebra is the same as for $M(H) = H^{\text{op}} \triangleright\!\blacktriangleleft H$. We have that $\theta$ is also an isomorphism of coalgebras, proving thus that $H^{\text{op}} \ {}^{\psi}\!\!\triangleright\!\blacktriangleleft H_\chi$ is a Hopf algebra. Here, its cross coproduct is explicitly :

$$\Delta(h \otimes a) = h_{(1)} \otimes h_{(2)(1)}Sh_{(2)(3)}h_{(3)(1)}\chi^{(1)}Sh_{(3)(4)}a_{(1)}\chi^{-(1)} \otimes h_{(2)(2)} \otimes h_{(3)(2)}\chi^{(2)}Sh_{(3)(3)}a_{(2)}\chi^{-(2)}$$

$$= h_{(1)} \otimes h_{(2)}\chi^{(1)}Sh_{(6)}a_{(1)}\chi^{-(1)} \otimes h_{(3)} \otimes h_{(4)}\chi^{(2)}Sh_{(5)}a_{(2)}\chi^{-(2)} \qquad (5)$$

We leave it to the reader to verify that

$$\Delta\theta(h \otimes a) = \theta \otimes \theta(h_{(1)} \otimes \chi^{(1)}a_{(1)}\chi^{-(1)} \otimes h_{(2)} \otimes \chi^{(2)}a_{(2)}\chi^{-(2)}) \qquad (6)$$

as required. $\square$

There are many examples of twisting of interest. for example, if $H_0$ is some quasitriangular Hopf algebra, then the quasitriangular structure $\mathcal{R} \in H_0 \otimes H_0$ can be viewed as a cocycle $\chi = \mathcal{R}_{23}$ where we view $\mathcal{R}$ in $H_0 \otimes H_0 \otimes H_0 \otimes H_0$ in the second and third factors. Then $H = H_0 \otimes H_0$ and $M_\chi(H)$ is a cocycle bicrossproduct with

$$\psi(h \otimes g) = h_{(1)} Sh_{(4)} \otimes g_{(1)} \mathcal{R}^{(1)} Sg_{(2)} \mathcal{R}^{-(1)} \otimes h_{(2)} \mathcal{R}^{(2)} Sh_{(3)} \mathcal{R}^{-(2)} \otimes 1 \quad (7)$$

Some physical applications of the hereabove result have been investigated in [4,5,6]. In particular, it is known that the q-Lorentz group $Uq(so(3, 1))$ is a twisting of the q-Euclidean group $Uq(so(4))$ (see [10]). Consequeltly, starting from the general form $M_\chi(H)$, the above $M\chi(Uq(so(3, 1)))$ provides a single Hopf algebra of the form $Uq(so(4))^{op\,\psi} \bowtie Uq(so(3, 1))$ unifying the Lorentzian and Euclidean ones. It is noted that the case $Uq(so(2,2))$ is exluded. This result may represent an interesting algebraic framework for the hypothesis of fluctuation of the signature of the spacetime metric at the Planck scale between Lorentzian and Euclidean signatures proposed in (4).

## 3. THE COCYCLE BICROSSPRODUCT $\overline{M}(H) = H^\psi \bowtie H$

Here we give now a related construction that was originally obtained as an example of the above in the quasitriangular case. However, the result may be proven directly without the quasitriangular structure and holds more generally.

**Proposition 3.1** *Let H be any Hopf algebra equipped with a bijective antipode. There exists a cocycle bicrossproduct $H^\psi \bowtie H$ where*

$a \triangleleft h = Sh_{(1)} a\ h_{(2)}$
$\beta(h) = Sh_{(1)} h_{(3)} \otimes h_{(2)}$
$\psi(h) = Sh_{(1)} h_{(3)} \otimes Sh_{(2)} h_{(4)}$

*In addition, $H^\psi \bowtie H \cong H \otimes H$, as Hopf algebras.*

**Proof** Clearly $H$ acts on $H$ by the right quantum adjoint action as usual, so we have an algebra $H \bowtie H$. Meanwhile, $\beta$ as stated is the usual right adjoint coaction viewed as a left cocycle action.

For this to work one may check that one needs the dual cocycle $\psi$ as stated. One may check that one then fulfills all the conditions for a cocycle bicrossproduct. Alternatively, we let $\theta$ be defined by

$$\theta : H \otimes H \xrightarrow{\sim} H^\psi {\triangleright\blacktriangleleft} H$$

$$\theta(h \otimes g) = h_{(1)} \otimes Sh_{(2)}g$$

and show that this gives the required isomorphism. For the algebra $H \mathbin{\triangleright\!\!\!<} H$ this is standard and the same as for $M(H)$ in [11]. Less trivial, the coproduct of $H^{\,\psi}{\triangleright\blacktriangleleft} H$ is :

$$\Delta h \otimes a = h_{(1)} \otimes Sh_{(2)(1)}h_{(2)(3)}\psi(h_{(3)})^{(1)}a_{(1)} \otimes h_{(2)(2)}\psi(h_{(3)})^{(2)}a_{(2)}$$

$$= h_{(1)} \otimes Sh_{(2)}h_{(5)}a_{(1)} \otimes h_{(3)} \otimes Sh_{(4)}h_{(6)}a_{(2)}$$

and we check :

$$\Delta\theta(h \otimes a) = h_{(1)(1)} \otimes Sh_{(1)(2)}h_{(1)(5)}Sh_{(2)(2)}a_{(1)} \otimes h_{(1)(3)} \otimes Sh_{(1)(4)}h_{(1)(6)}Sh_{(2)(1)}a_{(2)}$$

$$= h_{(1)} \otimes Sh_{(2)}a_{(1)} \otimes h_{(3)} \otimes Sh_{(4)}a_{(2)}$$

$$= (\theta \otimes \theta)((h_{(1)} \otimes a_{(1)}) \otimes (h_{(2)} \otimes a_{(2)}))$$

$$= (\theta \otimes \theta)\Delta_{H \otimes H}(h \otimes a)$$

as required. $\square$

Note that if $H$ is quasitriangular, then the quasitriangular structure $\mathcal{R}$ could be viewed as a cocycle and $H_\chi = H^{cop}$ whith the opposite coproduct. Then the construction of theorem 2.1 gives a bicrossproduct $(H^{op})^\psi {\triangleright\blacktriangleleft} H^{cop}$ which coincides with $\overline{M}(H^{cop})$ on noting that $H^{cop}$ isomorphic to $H^{op}$ via the Hopf algebra antipode.

## 4. CONCLUSION

Finally, we mention that we would expect the cocycle bicrossproducts of theorem 2.1 and theorem 3.1 corresponds to quasiassociative (or coquasi-Hopf algebras) by semidualising one of the factors. In general, one expects that a cocycle bicrossproduct $H \triangleright\blacktriangleleft A$ with $A$ finite dimensional should semidualise to a generalise double crossproduct $A^* \triangleright\triangleleft H$ as a dual quasi-Hopf algebra. This is along the same lines as without cocycles in [1] and was conjectured already in [3]. In particular, we expect a coquasi-Hopf algebra $A \underset{\chi}{\triangleright\triangleleft} H$ where $A$ is dually paired with $H$. Details will be developed elsewhere, but we note that a general result linking Hopf algebra extensions and quasi-Hopf algebras via monoïdal categories recently appeared in [12] could provide an intersting context for these constructions. Physically also, one can observe that Hopf algebra duality and the Drinfeld double at the Lie bialgebra level has been related to T-duality in string theory [13].

## REFERENCES

1. Majid S "Physics for algebraists : non-commutative and non-cocommutative Hopf algebras by a bicrossproduct construction"  *Journal of Algebra* **130**  (1990)

2. Majid S "More examples of bicrossproducts and double crossproducts Hopf algebras" *Isr. Jour. Math.* **72**  (1990)

3. Majid S and Schroers B  « q-Deformation and semidualisation in 3D quantum gravity »  *J. Phys A* **42** (2009)

4. Bogdanoff G "Fluctuations quantiques de la signature de la métrique à l'échelle de Planck" *Thèse de Doctorat Univ. de Bourgogne* (1999)

5. Bogdanoff G and  Bogdanoff I " Topological field theory of the initial singularity of spacetime" *Classical and  Quantum Gravity*  **21**  (2001)

6. Bogdanoff G and  Bogdanoff I. "Spacetime metric and KMS condition" *Annals of Physics* **295** (2002)

7. Drinfeld V.G  "On almost-cocommutative Hopf algebras"  *Leningrad Math. J.* **1**  (1990)

8. Gurevich D.I   Majid S "Braided groups of Hopf algebras obtained by twisting" *Pac. J. Math* **162** (1994)

9. Drinfeld  V.G  In A. Gleason Ed.  "Quantum groups" *Proceedings of the ICM*. Rhode Island (1987).

10. Majid  S. "q-Euclidean space and quantum Wick rotation by twisting". *J. Math. Phys.* **35**, 9 (1994)

11 . Majid S "Foundations of quantum groups" *Cambridge University Press*  (1995)

12. Schauenberg  P "Hopf algebra extensions and monoïdal categories" Personal communication of the author. *http : / / www. mathematik . uni-muchen.de / personen / schauenberg . html*  (2001)

13. Klimcik C and Severa P "Poisson Lie T-duality and loop groups of Drinfeld doubles" *Phys. Let. B* **372** (1996)